\documentclass[12pt]{amsart}
\numberwithin{equation}{section}

\usepackage{amssymb}
\def\ca{{\mathcal A}}

\def\cd{{\mathcal D}}

\def\cf{{\mathcal F}}

\def\ch{{\mathcal H}}

\def\cn{{\mathcal N}}

\def\cs{{\mathcal S}}

\def\cz{{\mathcal Z}}

\def\ga{{\mathfrak A}} \def\gpa{{\mathfrak a}}
\def\gb{{\mathfrak B}}

\def\gf{{\mathfrak F}}
 \def\gpg{{\mathfrak g}}

 \def\gps{{\mathfrak s}}
 \def\gpt{{\mathfrak t}}

\def\gz{{\mathfrak Z}}

\def\br{{\mathbb R}}
\def\bt{{\mathbb T}}

\def\bz{{\mathbb Z}}

\def\a{\alpha}
\def\b{\beta}
\def\g{\gamma} \def\G{\Gamma}
\def\d{\delta}  
\def\eeps{\epsilon}
\def\eps{\varepsilon}
 \def\L{\Lambda}

\def\m{\mu}

\def\n{\nu}
\def\r{\rho}
\def\s{\sigma} \def\S{\Sigma}
\def\t{\tau}
\def\f{\varphi}\def\ff{\phi} 
\def\th{\theta}  
\def\om{\omega} \def\Om{\Omega}
\def\z{\zeta}

\def\id{\hbox{id}}

\newtheorem{thm}{Theorem}[section]

\newtheorem{cor}[thm]{Corollary}
\newtheorem{prop}[thm]{Proposition}

\def\aut{\mathop{\rm Aut}}

\def\esssup{\mathop{\rm esssup}}

\def\sp{\mathop{\rm sp}}

\newcommand{\ty}[1]{\mathop{\rm {#1}}}
\def\di{\mathop{\rm d}\!}
\def\vol{\mathop{\rm vol}}

\newcommand{\nn}{\nonumber}

\begin{document}

\title[chemical potential for disordered systems]
{KMS states and the chemical potential for disordered systems}
\author{Francesco Fidaleo}
\address{Francesco Fidaleo,
Dipartimento di Matematica,
Universit\`{a} di Roma Tor Vergata, 
Via della Ricerca Scientifica 1, Roma 00133, Italy} \email{{\tt
fidaleo@mat.uniroma2.it}}
\date{\today}

\begin{abstract}
We extend the theory of the chemical potential associated to a compact 
separable gauge group to the case of disordered quantum systems. 
This is done in the natural framework of operator algebras. Among the 
other results, we show that the chemical potential does not depend on 
the disorder. The situation of the $n$--torus is treated in some 
detail. Indeed, provided that the zero--point is fixed independently 
on the disorder,
the chemical potential is intrinsically defined in 
terms of the direct integral decomposition of the
Connes--Radon--Nikodym cocycle associated to the KMS state
$\om$ and its trasforms $\om\circ\r$ by the localized automorphisms 
$\r$ of the observable algebra, carrying the abelian charges of the 
model under consideration. This description parallels the analogous 
one relative to the usual (i.e. non disordered) quantum models. 
\end{abstract}

\maketitle

\section{introduction}
\label{sec1}

In quantum physics, one often recovers the observable algebras by a 
principle of global gauge invariance. The reader is referred to 
\cite{DHR1, DHR2, DR1, DR2} and the reference cited therein. In order 
to investigate the termodynamical behavior of such physical models, 
the concept of chemical potential naturally arises. The algebraic 
description of the chemical potential is well understood, taking into 
account the principle of gauge invariance. Namely, suppose that we 
have $k$ species of particles (i.e. the chemical components). 
Provided that the 
field algebra $\cf$ has a natural local structure, one can 
consider on it, infinite volume limits of states 
arising from the Gibbs grand canonical ensamble relative to fixed 
inverse temperature $\b$ and the $k$--parametric chemical potential 
$\m=(\m_{1},\dots,\m_{k})$, one for 
each species of particles. It is seen that each of these 
states $\f_{\b,\m}$ satisfies the Kubo--Martin--Schwinger 
(KMS for short) boundary condition for a one parameter subgroup of 
the $(k+1)$ dimensional Lie group $\br\times G$, the gauge 
group $G$ being isomorphic to the $k$--dimensional 
torus in this situation.\footnote{Here, in order to simplify matter, we are supposing 
that the asymmetry subgroup of a KMS state on the field algebra is 
trivial, see Section \ref{secesds} for the 
asymmetry subgroup.} If we denote by $H$ and 
$\{N_{i}\}$ the infinitesimal generators of the time translations and 
the gauge transformations respectively, the one parameter subgroup 
mentioned above has the form
$$
X=H-\sum_{i}\m_{i}N_{i}\,.
$$

As the observables are gauge invariant, the restriction 
$\om_{\b,\m}$
of any $\f_{\b,\m}$ as above to the observable algebra 
$\ca=\cf^{G}$, satisfies the KMS condition at the same inverse 
temperature $\b$ for the time evolution. In general, states corresponding to different 
values of the chemical potentials $\{\m_{i}\}$, give rise to 
different KMS states, when restricted to the observable algebra.

In quantum physics, all the physical content of the model is encoded 
directly in the observable algebra. On the one hand, the chemical 
potential has a physical meaning, even if it naturally arises by the 
use of the field algebra. On the other hand, the KMS boundary 
condition does not refer to the local structure of the algebra of the 
observable. 

In the seminal paper \cite{AKTH}, all these questions are 
explained in detail, see also \cite{AK, KT1, KT2} for strictly connected questions. 
Namely, if $\om$ is an extremal (or, more generally 
a weakly clustering) KMS state on $\ca$,
it is shown that any weakly clustering extension to all of $\cf$ is a 
KMS state relative to a new evolution modified by a one parameter 
subgroup of the gauge group. In the simple situation described above, 
this one parameter subgroup uniquely determines the values of the 
chemical potentials. 

Unfortunately, the results of \cite{AKTH} are not directly applicable 
to disordered models, the last including the very interesting examples 
of the spin glasses. The
equilibrium statistical mechanics of models arising from spin 
classical glasses has been
intensively studied, in order to understand the complex behaviour of 
the set of its temperature states. We refer the reader, for example, 
to \cite{B1, EH, MPV, Ne, NeS1} and the literature cited
therein. Some attempts to understand the structure of the set of the 
KMS states of quantum disordered systems is made in \cite{A1, B, BF, K} 
by using the standard techniques of operator algebras. In the present 
paper, we follow the last strategy.

Namely, in order to achieve the disorder, it is natural to set for the 
observable algebra, $\ga:=\ca\otimes L^{\infty}(X,\n)$, $(X,\n)$ 
being the sample space for the coupling ``constants'' of the system. 

As it is noted in \cite{BF}, most of the interesting states of a 
disordered system are states $\om$ whose centre $\gz_{\pi_{\om}}$ of 
the GNS representation $\pi_{\om}$ contains an Abelian algebra which 
is isomorphic to $L^{\infty}(X,\n)$. In addition, the KMS states of 
interest 
for which $\gz_{\pi_{\om}}$ is precisely $L^{\infty}(X,\n)$, can be 
interpreted as the ``pure termodynamical phases'' in the case of 
disordered systems. Furthermore, the phenomenon of the ``weak 
Gibbsianess'' naturally appear in this case. Indeed, if one consider 
infinite volume limits of finite volume Gibbs states, one obtain 
states on $\ca\otimes L^{\infty}(X,\n)$ satisfying the following 
properties. Its marginal 
distribution of the couplings is the given probability measure 
describing the disorder, whereas the conditional distribution of the 
standard observable variables (the ``spin'' variables,), given the couplings, 
is some 
infinite--volume Gibbs state 
almost surely. Such a field of infinite--volume  
Gibbs states satisfies an equivariant 
property (see \eqref{eqiv}). Namely, 
is it gives a (quantum version 
of a) Aizenman--Wehr metastate (see \cite{AW}), after direct integral decomposition, 
the last one assuming the meaning  
of the quantum counterpart of 
the classical procedure of conditioning w.r.t. 
the disorder variables. For recent 
results on metastates, we refer the reader to \cite{Ne, NeS1} 
and the references cited therein. It can happen that the state so 
obtained is not necessarily jointly Gibbsian relatively to the 
standard variables and the couplings, see \cite{EMSS, Ku} for some 
pivotal classical examples. 

The systematic investigation of the 
difference between Gibbsianess and weak Gibbsianess (i.e. states 
which arise from infinite volume limit of finite volume Gibbs states 
but are not jointly Gibbsian) started in the paper \cite{K}
in the setting of operator algebras. In this paper, both equilibrium 
conditions are connected with some natural variational principles. 
As it is explained in Section 6 of 
\cite{BF}, the KMS boundary condition for the algebra generated by 
spin variables and disorder variables    
seems to describe the weak 
Gibbsianess in quantum case. 
Moreover, it does not refer to the local structure of the observable algebra. 
For the r\^ole  played by Gibbsianess and weak Gibbsianess 
in the description of the termodynamical behavior of 
a disordered model, we refer the reader to the above mentioned papers.

In the present paper we extend the algebraic description of the 
chemical potential to disordered systems without referring to the 
difference between Gibbsianess and weak Gibbsianess. We take 
advantage twice by the paper \cite{AKTH}. First, we follow its 
plan. Second, we use the results of this paper in order to describe 
the occurrence of the chemical potential for disordered models.

The present paper is organized as follows. After a preliminary section, 
in Section \ref{secerg} we 
investigate some useful ergodic propoerties of states of interest of 
disordered systems. Section \ref{secesds} is devoted to the occurrence 
of the chemical potential. Starting from a state $\f$ on the field 
algebra $\gf=\cf\otimes L^{\infty}(X,\n)$, normal when restricted to 
the subalgebra $L^{\infty}(X,\n)$, which is weakly clustering 
with respect to the spatial translations, we show that the 
stabilizer, as well as the asymmetry subgroup coincide almost surely 
with the corresponding objects relative to the states $\f_{\xi}$. 
Here, the 
measurable equivariant field $\{\f_{\xi}\}_{\xi\in X}\subset\cf$ 
provides the 
direct integral decomposition of $\f$. Then, we show that for any weakly 
clustering state $\f$ on $\gf$ whose restriction to $\ga$ is KMS,
there exists a modification of the time 
evolution by a suitable one parameter group of the gauge group, the same for 
each $\f_{\xi}$, such that $\f_{\xi}$ is KMS with respect to this 
modified evolution almost surely. This is the content of Theorem 
\ref{gaga2} which is the natural generalization to our situation 
of Theorem II.4 of \cite{AKTH}. Section \ref{secint} is devoted to an 
intrinsic description of the chemical potential directly in terms of 
objects related to the algebra of observables. The case of the unit 
circle is treated in some detail, the case of the $n$--torus being 
quite similar. 

Provided that its zero--point is fixed independently 
on the disorder,
the chemical potential is intrinsically defined in 
terms of the Connes--Radon--Nikodym cocycle associated to the KMS state
$\om\in\cs(\ga)$ under consideration,
and its trasforms $\om\circ\r$ by the localized automorphisms 
$\r$ of the observable algebra $\ga$, carrying the abelian charges of the 
model. Indeed, under suitable conditions, it is proven that 
$\om\circ\r$ is equivalent to $\om$ also in our situation. 
Furthermore, the chemical potential is connected, and is independent 
almost surely on the disorder, with the Connes--Radon--Nikodym cocycle
$\big(D(\om_{\xi}\circ\r_{\xi}):D\om_{\xi}\big)$ of 
$\om_{\xi}\circ\r_{\xi}$ relative to $\om_{\xi}$, see Formula \eqref{crncoc1}. 
Here, $\{\om_{\xi}\}_{\xi\in X}$ provides the direct integral decomposition of 
(the normal extension of) $\om$, and the measurable field of normal automorphisms 
$\{\r_{\xi}\}_{\xi\in X}$ give rise the normal extension of $\r$ to all of the 
$\pi_{\om_{\xi}}(\ca)''$ which exist by Proposition \ref{esdi}.

\section{preliminaries}
\label{sec2}

We start by recalling the definition of the KMS boundary condition.
A state $\ff$ on the $C^{*}$--algebra $\gb$ satisfies the KMS
boundary condition at inverse temperature $\b$ which we suppose to be 
always different from zero, 
w.r.t the group of
automorphisms $\{\t_{t}\}_{t\in\br}$ if
\begin{itemize}
\item[(i)] $t\mapsto\ff(A\t_{t}(B))$ is a continuous function for
every $A,B\in\gb$,
\item[(ii)]
$\int\ff(A\t_{t}(B))f(t)\di t=\int\ff(\t_{t}(B)A)f(t+i\b)\di t$
whenever $f\in\widehat{\cd}$, $\cd$ being the 
space made of all infinitely 
often differentiable compactly supported functions in $\br$.
\end{itemize}

For the equivalent characterizations of the KMS boundary
condition, the main results about KMS states, and finally the
connections with Tomita theory of von Neumann algebras, see
e. g. \cite{BR2, St} and the references cited therein.

It is well--known that the cyclic vector $\Om_{\ff}$ of the GNS 
representation $\pi_{\ff}$ is also separating for $\pi_{\ff}(\gb)''$. 
Denote with an abuse of notation, $\s^{\ff}$ its modular group. 

According to this definition of KMS boundary condition, we have
\begin{equation}
\label{modgns}
\s^{\ff}_{t}\circ\pi_{\ff}=\pi_{\ff}\circ\t_{-\b t}\,.
\end{equation}

Our set--up is a separable $C^{*}$--algebra $\ca$ with an identity $I$,
describing the 
physical observables.\footnote{In order to avoid technical 
complications, in quantum field theory the local algebras of 
observables are enlarged by taking the weak operator closure in the 
vacuum representation. Namely, the local algebras of observables are 
typically von Neumann algebras with separable predual, the former being non 
separable $C^{*}$--algebras. This does not affect the substance of the 
theory, see the comments in Section \ref{secint}.} 
We suppose that $\ca$ is obtained 
as the fixed--point algebra $\ca=\cf^{G}$ under a pointwise--norm continuous 
action 
$$
\g:g\in G\mapsto\g_{g}\in\aut(\cf)
$$
of a compact second countable group $G$ (the {\it gauge group}) on 
another separable $C^{*}$--algebra $\cf$ (the {\it field algebra}). 
This is a typical situation appearing in quantum field theory, when 
the charges present in the model are described in terms of
a principle of (global) 
gauge invariance, see e.g. \cite{DHR1, DHR2, DR1, DR2}. The present 
description can be 
applied also to 
nontrivial models where the local algebras of 
observables are full matrix algebras, see e.g. \cite{NS}
for a possible example along this line. 

We suppose that the group $\{\a_{x}\}_{x\in\bz^{d}}$ of spatial translations 
acts on $\cf$.
We consider also a standard measure space $(X,\n)$ 
based on a compact separable space $X$, and a Borel probability 
measure $\n$. The group $\bz^{d}$ of the spatial translations
is supposed to act on the probability space $(X,\n)$
by measure preserving ergodic transformations
$\{T_{x}\}_{x\in\bz^{d}}$.

A one parameter random group of automorphisms 
\begin{equation*}
%\label{5} 
(t,\xi)\in\br\times X\mapsto\t^{\xi}_{t}\in\aut(\cf)
\end{equation*}
is acting on $\cf$. It is supposed to be strongly continuous in the 
time variable for each 
fixed $\xi\in X$, and jointly strongly measurable. Consider, for 
$A\in\cf$, the 
strongly measurable function $f_{A,t}(\xi):=\t^{\xi}_{t}(A)$. We get
\begin{equation*}
\|f_{A,t}\|_{L^{\infty}(X,\n;\cf)}\equiv
\esssup_{\xi\in X}\|\t_{t}^{\xi}(A)\|_{\cf}
=\|A\|_{\cf}\,,
\end{equation*}
where the last equality follows as $\t^{\xi}_{t}$ is isometric.
We assume further 
that $\t$ acts locally. Namely, if $A$ is an element of $\cf$, 
then the function 
$f_{A,t}\in L^{\infty}(X,\n;\cf)$ 
belongs to the $C^{*}$--subalgebra 
$\cf\otimes L^{\infty}(X,\n)$, where the above $C^{*}$--tensor product 
is uniquely determined as
any commutative $C^{*}$--algebra is nuclear.

We assume the following commutation rules
\begin{align}
\label{4} 
\t^{T_{x}\xi}_{t}\a_{x}=&\a_{x}\t^{\xi}_{t}\nn\\
\a_{x}\g_{g}=&\g_{g}\a_{x}\\
\t^{\xi}_{t}\g_{g}=&\g_{g}\t^{\xi}_{t}\nn
\end{align}
for each $x\in\bz^{d}$, $\xi\in X$, $t\in\br$, and $g\in G$.

By \eqref{4}, it is immediate to show that $\a_{x}$ and $\t^{\xi}_{t}$ 
leave globally stable $\ca$. Namely, $\bz^{d}$, $\br$ act on $\ca$ as groups 
of automorphisms or random automorphisms, respectively. 

Finally, we address also the situation when Fermion operators are 
present in $\cf$. Namely, there exists an automorphism $\s$ of $\cf$ 
commuting with the all the gauge transformations, the spatial 
translations and the random time evolution, such that 
$\s^{2}=e$.\footnote{In most of the interesting physical 
situations, $\s\in\cz(G)$, $G$ being the gauge group, see e.g. \cite{DHR1, DR2}.
The situation without Fermion operators corresponds to $\s=e$.} We put 
\begin{equation}
\label{fesi}
\cf_{+}:=\frac{1}{2}(e+\s)(\cf)\,,\quad
\cf_{-}:=\frac{1}{2}(e-\s)(\cf)\,.
\end{equation}

The disordered system under consideration is described by
$$
\gf:=\cf\otimes L^{\infty}(X,\n)\,.
$$

Notice that, by identifying $\gf$ with a closed subspace of
$L^{\infty}(X,\n;\cf)$, each element $A\in\gf$ is uniquely
represented by a measurable essentially bounded function
$\xi\mapsto A(\xi)$ with values in $\cf$.

The group $\bz^{d}$ of all the space translations is naturally acting on
the $C^{*}$--algebra $\gf$ as 
\begin{equation*}
\gpa_{x}(A)(\xi):=\a_{x}(A(T_{-x}\xi))\,.
\end{equation*}
Further, define on $\gf$, 
\begin{align}
\label{aztemp1}
&\gpt_{t}(A)(\xi):=\t_{t}^{\xi}(A(\xi))\,,\nn\\
&\gpg:=\g\otimes\id_{L^{\infty}(X,\n)}\,,\\
&\gps:=\s\otimes\id_{L^{\infty}(X,\n)}\,.\nn
\end{align}

It is straightforward to verify that $\{\gpa_{x}\}_{x\in\bz^{d}}$,
$\{\gpt_{t}\}_{t\in\br}$ and $\{\gpg_{g}\}_{g\in G}$ define actions 
of $\bz^{d}$, $\br$ and $G$ on $\gf$ which are mutually commuting, 
and commute also with the parity automorphism $\gps$. The subspaces 
$\gf_{+}$ and $\gf_{-}$ are defined as in \eqref{fesi}.
Furthermore, taking into account \eqref{4} and the definition \eqref{aztemp1}
of the action of the gauge group on the disordered field algebra $\gf$, 
$\gpa_{x}$ and $\gpt_{t}$ leave globally stable 
the disordered observable algebra $\ga$. Namely, 
$\{\gpa_{x}\}_{x\in\bz^{d}}$ and $\{\gpt_{t}\}_{t\in\br}$ define by 
restriction, 
mutually commuting actions of $\bz^{d}$ and $\br$ on $\ga$, 
respectively. 

In order to study a class of states of interest for disordered 
systems, we start with  $*$--weak measurable fields of states
$$
\xi\in X\mapsto\f_{\xi}\in\cs(\cf)\,.
$$

We suppose that the field $\{\f_{\xi}\}_{\xi\in X}$ fulfils 
almost surely, the equivariance condition
\begin{equation}
\label{eqiv}
\f_{\xi}\circ\a_{x}=\f_{T_{-x}\xi}
\end{equation}
w.r.t. the spatial translations, simultaneously.
 
A state $\f$ on $\gf$ is naturally defined as follows:
\begin{equation}
\label{ddf}
\f(A)=\int_{X}\f_{\xi}(A(\xi))\n(\di\xi)\,,
\qquad A\in\gf\,.
\end{equation}

It is immediate to verify that $\f$ defined as above is invariant 
w.r.t. the space translations $\gpa_{x}$. Moreover, 
$\f\lceil_{I\otimes L^{\infty}(X,\n)}$ is a normal state.

Equally well, one can start with a  $\gpa$--invariant state $\f$ on $\gf$, which 
is normal when restricted to $I\otimes L^{\infty}(X,\n)$. Then, we can 
recover a $*$--weak measurable field $\{\f_{\xi}\}_{\xi\in X}\subset\cf$
fulfilling \eqref{eqiv}. Such a measurable fields provides the direct 
integral decomposition of $\f$ as in \eqref{ddf}, see \cite{BF}, Theorem 4.1. Similar
considerations can be applied to the observable algebras 
$\ga$ as well. In the 
sequel, we denote by $\cs_{0}(\ga)$, $\cs_{0}(\gf)$ the convex closed 
subset of states on $\ga$, $\gf$ respectively, fulfilling the 
properties listed above.

\section{ergodic properties of states of disordered systems}
\label{secerg}

In this section we study some useful ergodic properties of states in 
$\cs_{0}(\gf)$ or $\cs_{0}(\ga)$. We restrict ourselves to the field 
algebra, the other case being similar.  

Let $C,D\in\gf$, and $A,B\in\gf_{+}\bigcup\gf_{-}$. Put $\eeps_{A,B}=-1$ 
if $A,B\in\gf_{-}$ and $\eeps_{A,B}=1$ in the three remaining 
possibilities.
We say that the state $\f$ is 
{\it asymptotically Abelian} w.r.t. $\gpa$ if
\begin{equation}
\label{aaf}
\lim_{|x|\to+\infty}
\f\left(C\big(\gpa_{x}(A)B-\eeps_{A,B}B\gpa_{x}(A)\big)D\right)=0\,,
\end{equation}

The state $\f$ is 
{\it weakly clustering} w.r.t. $\gpa$ if
\begin{equation}
\label{wcf}
\lim_{N}\frac{1}{|\L_{N}|}\sum_{x\in\L_{N}}
\f(A\gpa_{x}(B))=\f(A)\f(B)\,,
\end{equation}
$\L_{N}$ being the box with a vertex sited in the 
origin, containing $N^{d}$ points with positive 
coordinates.\footnote{Taking into account natural 
applications to continuous disordered systems, 
we can consider cases when $\br^{d}$ is 
acting as the group of spatial translations. We use 
in \eqref{wcf} the natural 
modification $M$ of the Cesaro mean given on bounded measurable functions, by
$$
M(f):=\lim_{D\to+\infty}\frac{1}{\vol(\L_{D})}\int_{\L_{D}}f(x)\di^{d}x\,,
$$
$\L_{D}$ being a box with edges of lenght $D$. 
Most of the forthcoming analysis can 
be applied in this situation as well.}

Notice that, lots of interesting states are naturally asymptotically Abelian 
w.r.t. the spatial translations, see e.g. \cite{BF}, Proposition 
2.3, see also \cite{L}.

We report the following result for the 
sake of completeness.
\begin{prop}
\label{folk}
Suppose that $\f\in\cs(\gf)$ is a
$\gpa$--invariant asymptotically Abelian state. Then the following 
assertions are equivalent.
\begin{itemize}
\item[(i)] $\f$ is  $\gpa$--weakly clustering,
\item[(ii)] $\f$ is  $\gpa$--ergodic.
\end{itemize}
\end{prop}
\begin{proof}
It is a well--known fact that (i) always implies (ii). The reverse implication 
follows as in 
Proposition 5.4.23 of \cite{BR2}, the last working also under 
the weaker condition \eqref{aaf}.
\end{proof}

One sees that an 
asymptotically abelian invariant state is automatically 
$\gps$--invariant, that is it is an even state. 

The weak clustering property for states in $\cs_{0}(\gf)$ can be 
translated as a property of the corresponding equivariant fields of 
states on $\cf$. Namely, Let $\f\in\cs_{0}(\gf)$, and 
$\{\f_{\xi}\}_{\xi\in X}\subset\cs(\cf)$ the corresponding $\a$--equivariant 
measurable field of states. Then the results listed below hold true.
\begin{prop}
The following assertions are equivalent.
\begin{itemize}
\item[(i)] $\f$ is $\gpa$--weakly clustering,
\item[(ii)] we have for each $A\in\cf$, $B\in\gf$,
\begin{equation}
\label{wcf1}
\lim_{N}\frac{1}{|\L_{N}|}\sum_{x\in\L_{N}}
\f_{\xi}(A\a_{x}(B(T_{-x}\xi)))=\f_{\xi}(A)\f(B)\,,
\end{equation}
in the $*$--weak  topology of $L^{1}(X,\n)$.
\end{itemize}
\end{prop}
\begin{proof}
(i) $\Rightarrow $(ii) We compute for $f\in L^{\infty}(X,\n)$, $A\in\cf$,
$B\in\gf$,
\begin{align*}
&\int_{X}f(\xi)(\f_{\xi}(A)\f(B))\equiv\f(A\otimes f)\f(B)\\
=&\lim_{N}\frac{1}{|\L_{N}|}\sum_{x\in\L_{N}}\f(A\otimes f\gpa_{x}(B))\\
=&\lim_{N}\int_{X}f(\xi)
\bigg(\frac{1}{|\L_{N}|}\sum_{x\in\L_{N}}
\f_{\xi}(A\a_{x}(B(T_{-x}\xi)))\bigg)
\end{align*}
which is the assertion.

(ii) $\Rightarrow $(i) By a standard density argument, we reduce the 
situation to element of the form $A\otimes f,B\in\gf$, with
$f\in L^{\infty}(X,\n)$, $A\in\cf$. We get
\begin{align*}
&\lim_{N}\frac{1}{|\L_{N}|}\sum_{x\in\L_{N}}\f(A\otimes f\gpa_{x}(B))\\	
=&\lim_{N}\int_{X}f(\xi)
\bigg(\frac{1}{|\L_{N}|}\sum_{x\in\L_{N}}
\f_{\xi}(A\a_{x}(B(T_{-x}\xi)))\bigg)\\
=&\int_{X}f(\xi)(\f_{\xi}(A)\f(B))\equiv\f(A\otimes f)\f(B)
\end{align*}
and we are done.
\end{proof}

A sufficient condition for the weak clustering property of $\f$ is the 
pointwise--clustering property for the corresponding equivariance field
$\{\f_{\xi}\}_{\xi\in X}$ of states, almost surely. 
\begin{prop}
\label{ces1}
Suppose that, for each $A,B\in\cf$ 
\begin{equation}
\label{ces}
\lim_{|x|\to+\infty}\f_{\xi}(A\a_{x}(B))=\f_{\xi}(A)\f_{T_{-x}\xi}(B)
\end{equation}
almost surely. Then the state $\f$ given by \eqref{ddf} is weakly 
clustering.
\end{prop}
\begin{proof} By a standard density argument, we can reduce the 
situation to a measurable set $F\subset X$ of full measure such that
\eqref{ces} is satisfied simultaneously for each $A,B\in\cf$. Define, 
for fixed $A,B\in\cf$, $g\in L^{\infty}(X,\n)$ and $\xi\in F$,
\begin{align*}
\d_{1}(N):=&\frac{1}{|\L_{N}|}\sum_{x\in\L_{N}}g(T_{-x}\xi)
\bigg(\f_{\xi}(A\a_{x}(B))-\f_{\xi}(A)\f_{\xi}(\a_{x}(B))\bigg)\\
\d_{2}(N):=&\frac{1}{|\L_{N}|}\sum_{x\in\L_{N}}g(T_{-x}\xi)
\f_{T_{-x}\xi}(B)-\f(B\otimes g)\,.
\end{align*}
We have that $\d_{1}(N)\longrightarrow 0$ by hypotesis, and
$\d_{2}(N)\longrightarrow 0$ by the Individual Ergodic Theorem. Then we 
conclude that, for fixed $A,B\in\cf$, $g\in L^{\infty}(X,\n)$,
\begin{align*}
\frac{1}{|\L_{N}|}&\sum_{x\in\L_{N}}g(T_{-x}\xi)
\f_{\xi}(A\a_{x}(B))-\f_{\xi}(A)\f(B\otimes g)\\
&\equiv\d_{1}(N)-\f_{\xi}(A)\d_{2}(N)\longrightarrow 0
\end{align*}
pointwise on $F$, as $N\longrightarrow+\infty$.
Let now $f\in L^{\infty}(X,\n)$. We have
\begin{align*}
&\frac{1}{|\L_{N}|}\sum_{x\in\L_{N}}\f(A\otimes f\gpa_{x}(B\otimes g))-
\f(A\otimes f)\f(B\otimes g)\\
=&\int_{X}f(\xi)\bigg(\frac{1}{|\L_{N}|}\sum_{x\in\L_{N}}g(T_{-x}\xi)
\f_{\xi}(A\a_{x}(B))-\f_{\xi}(A)\f(B\otimes g)\bigg)\n(\di\xi)
\end{align*}
which goes to $0$ by the Lebesgue Dominated Convergence Theorem.
\end{proof}  

Finally, we point out the fact that there exist examples of weakly clustering 
states $\f$ satisfying 
the properties listed above without assuming any factoriality 
condition about $\f$, the last being not natural in the setting of 
disordered systems, see Section 4 and Section 5 of \cite{BF}. Indeed, 
it is enough to consider any quasi--local algebra $\cf$ admitting a state $\om$ 
which is strongly clustering (i.e. mixing) w.r.t. the space 
translations. Take any probability space $(X,\n)$ on which the space 
translations act ergodically. Define $\f_{\xi}=\om$, the constant field.
By Proposition \ref{ces1}, the state given by 
$\f(A):=\int_{X}\om(A(\xi))\n(\di\xi)$ satisfies all the required 
properties. In addition, it follows By Proposition 2.3 of \cite{BF},
that the states
on $\gf$ constructed as above from states $\om$ arising from quasi 
local algebras describing spin models on $\bz^{d}$ are 
$\gpa$--asymptotically Abelian as well.

We speak, without any further mention, and if it is not otherwise 
specified, about asymptotic Abelianity, weak 
clustering, or ergodicity for states, if 
they satisfy these properties w.r.t. the spatial translations.

\section{extension of states for disordered systems}
\label{secesds}

In this section we prove, following the line of \cite{AKTH}, the 
following result. For
each weakly 
clustering state $\f\in\cs_{0}(\gf)$ whose restriction to $\ga$ 
satisfies the KMS boundary condition,
there exists a modification of the time 
evolution by a suitable one parameter group of the gauge group, the same for 
each $\f_{\xi}$, such that $\f_{\xi}$ is KMS w.r.t. this 
modified evolution almost surely. This is done after showing that the 
stabilizer $G_{\f}\subset G$, as well as the asymmetry subgroup 
$N_{\f}\subset G$ coincide almost surely 
with the corresponding objects relative to the states $\f_{\xi}$.

Let $\f,\psi\in\cs_{0}(\gf)$, 
and $\{\f_{\xi}\}_{\xi\in X}$, $\{\psi_{\xi}\}_{\xi\in X}$ be the 
corresponding $\a$--equivariant measurable fields of states on $\cf$.
\begin{prop}
\label{gaga}
If $\f,\psi\in\cs_{0}(\gf)$ are weakly clustering states whose restrictions to 
$\ga$ are equal, then there exist $g\in G$ and a measurable subset $F$ of 
full measure, such that $\xi\in F$ implies
$$
\psi_{\xi}(A)=\f_{\xi}(\g_{g}(A)),
$$
simultaneously for every $A\in\cf$.
\end{prop}
\begin{proof}
By Theorem II.1 of \cite{AKTH}, there exists $g\in G$ such that 
$\psi=\f\circ\gpg_{g}$. We compute for each $f\in L^{\infty}(X,\n)$ and 
$A\in\cf$,
$$
\int_{X}\n(\di\xi)f(\xi)\big(\psi_{\xi}(A)-\f_{\xi}(\g_{g}(A))\big)=0\,.
$$
This means that for any fixed $A\in\cf$ there exists a measurable set 
$F_{A}$ of full measure such that 
$$
\psi_{\xi}(A)=\f_{\xi}(\g_{g}(A))
$$
on $F_{A}$. Choose a dense countable subset $\cf_{0}\subset\cf$. Put
${\displaystyle F:=\bigcap_{A\in\cf_{0}}F_{A}}$. $F$ is a measurable 
subset of $X$ of full measure. For each $A\in\cf$ choose a sequence 
$\{A_{n}\}$ in $\cf_{0}$ converging to $A$. We have for $\xi\in F$,
$$
\psi_{\xi}(A)=\lim_{n}\psi_{\xi}(A_{n})=\lim_{n}\f_{\xi}(\g_{g}(A_{n}))
=\f_{\xi}(\g_{g}(A))\,.
$$
\end{proof}
\begin{prop}
\label{wce}
Suppose that each state in $\cs_{0}(\gf)$ is asymptotically Abelian. Then
any weakly clustering state $\om\in\cs_{0}(\ga)$ extend to a weakly 
clustering state $\f\in\cs_{0}(\gf)$.
\end{prop}
\begin{proof}
Let $\psi$ be any extension of $\om$. Then it is normal when 
restricted to $L^{\infty}(X,\n)$. It could be not 
$\gpa$--invariant.\footnote{Notice that, for the measurable 
field $\{\psi_{\xi}\}_{\xi\in X}$ of positive forms giving the 
decomposition of $\psi$,
$\psi_{\xi}(I)=1$ almost everywhere. So we have a decomposition of 
$\psi$ into a measurable, not necessarily equivariant, field of states.}
Let $m$ be any invariant mean on $\bz^{d}$ which exists as $\bz^{d}$ is 
amenable. Then $m\big(\{\psi(\gpa_{x}(A)\}\big)$ defines an invariant 
extension of $\om$, that is the compact convex subset of $\cs_{0}(\gf)$
consisting of all the $\gpa$--invariant 
extensions of $\om$ is nonvoid. Take any estremal element $\f$ in 
such compact convex set. As $\om$ is weakly clustering, it is extremal 
in the set of all $\gpa$--invariant states (i.e. $\gpa$--ergodic). 
As $\f$ is an estremal 
extension of $\om$, we conclude that $\f$ is itself $\gpa$--ergodic. 
By Proposition \ref{folk}, $\f$ is also weakly clustering under our 
assumptions.
\end{proof}

Notice that, by the results contained in \cite{BF}, Proposition 2.3 
of \cite{BF}, and the considerations in Section \ref{secint}, 
there are disordered models satisfying the 
assumptions of Proposition \ref{wce}.

We define in the usual way, the stabilizer of a state $\f\in\cs(\gf)$ as
$$
G_{\f}:=\big\{g\in G\,\big|\,\f\circ\gpg_{g}=\f\big\}\,.
$$
The stabilizers $G_{\f_{\xi}}$ of the $\f_{\xi}\in\cs(\cf)$ are defined 
analogously as
$$
G_{\f_{\xi}}:=\big\{g\in G\,\big|\,\f_{\xi}\circ\g_{g}=\f_{\xi}\big\}\,.
$$

The normalizer $\cn(H)$ and the centralizer $\cz(H)$ of a subgroup 
$H\subset G$ are defined in the usual way as
\begin{align*}
&\cn(H):=\big\{g\in G\,\big|\,gHg^{-1}=H\big\}\,,\\
&\cz(H):=\big\{g\in G\,\big|\,gh=hg\,,h\in H\big\}\,.
\end{align*}

We show that the stabilizers $G_{\f_{\xi}}$ of the $\f_{\xi}$ are almost surely 
independent on the disorder, and coincide with the 
stabilizer $G_{\f}$ of $\f$ almost everywhere. 
\begin{thm}
\label{stbb}
Let $\f\in\cs_{0}(\gf)$. Then there exists a measurable
set $F\subset X$ of full measure 
such that $\xi\in F$ implies $G_{\f_{\xi}}=G_{\f}$.
\end{thm}
\begin{proof}
Without loss of generality, we can suppose that \eqref{eqiv} holds 
true everywhere on $X$.\footnote{Let $X_{0}\subset X$ be the 
measurable set of full measure such that \eqref{eqiv} is 
simultaneously satisfied for $x\in\bz^{d}$. We reduce the situation to 
the measurable invariant set ${\displaystyle\bigcap_{x\in\bz^{d}}T_{x}X_{0}}$
of full measure.}
Choose a countable dense subset $\{h\}$ of $G$ which always exists by 
separability. Define 
\begin{align*}
V_{h,n}:=&\big\{g\in G\,\big|\,{\rm dist}(g,h)\leq 1/n\big\}\,,\\
X_{h,n}:=&\big\{\xi\in X\,\big|\,G_{\f_{\xi}}\cap V_{h,n}
\neq\emptyset\big\}\,,
\end{align*}
where ``dist'' is any metric on $G$ generating its topology.
Put
\begin{align*}
f_{A}(\xi,g):=&\f_{\xi}(\g_{g}(A))-\f_{\xi}(A)\,,\\
\G:=&\bigcap_{A\in\cf}f_{A}^{-1}(\{0\})\,.
\end{align*}

By separability, we can reduce the last intersection to a countably 
dense set of $\cf$, that is $\G$ is a measurable subset of $X\times G$.
Furthermore, it is immediate to check that
$$
X_{h,n}=P_{X}\big((X\times V_{h,n})\cap\G\big)\,,
$$
$P_{X}$ being the projection w.r.t. the first variable. 
This means that the sets $X_{h,n}$ are measurable. Taking into 
account also \eqref{eqiv}, they are also invariant under the action 
of the spatial translations. So, by ergodicity, we have that 
$\n(X_{h,n})$ is $0$ or $1$. 

Let $S_{n}\subset G$ be the set of the $h$ such that $\n(X_{h,n})=1$. 
Define
\begin{align*}
G_{1}:=&\bigcap_{n>0}\bigg\{\bigcup_{h\in S_{n}}V_{h,n}\bigg\}\,,\\
F_{1}:=&\bigcap_{n>0}\bigg\{\bigcap_{h\in S_{n}}X_{h,n}\bigg\}\,.
\end{align*}
Notice that $F_{1}$ is a measurable set of full measure.

We get that $\xi\in F_{1}$ implies that $G_{\f_{\xi}}=G_{1}$. Indeed, 
if $g\in G_{1}$, ${\displaystyle g\in\bigcup_{h\in S_{n}}V_{h,n}}$ 
for each $n$, which means for each $n$, 
${\rm dist}\big(G_{\f_{\xi}},g\big)<1/n$ whenever $\xi\in F_{1}$. Hence,
$\xi\in F_{1}$ implies $G_{1}\subset G_{\f_{\xi}}$. Conversely, if 
$g\in G_{\f_{\xi}}$ and $\xi\in F_{1}$, then 
${\displaystyle g\in\bigcup_{h\in S_{n}}V_{h,n}}$ 
for every $n$, which means $G_{\f_{\xi}}\subset G_{1}$. 

Let now $\cf_{0}\subset\cf$, $G_{0}\subset G_{\f}$ be dense countable 
subsets, then 
we obtain for each $f\in L^{\infty}(X,\n)$, $A\in\cf_{0}$ and
$g\in G_{0}$,
$$
\int_{X}\n(\di\xi)f(\xi)\big(\f_{\xi}(\g_{g}(A))-\f_{\xi}(A)\big)=0
$$
which means that one can choose a measurable set $F_{0}$ of full 
measure such that $\xi\in F_{0}$, $A\in\cf_{0}$ and $g\in G_{0}$ 
imply 
$$
\f_{\xi}(\g_{g}(A))=\f_{\xi}(A)\,,
$$
see the proof of Proposition \ref{gaga}.

Let now $g\in G_{\f}$ be fixed. Choose a sequence $\{g_{n}\}\subset G_{0}$
converging to $g$. We have for each $\xi\in F_{0}$ and $A\in\cf_{0}$,
$$
\f_{\xi}(\g_{g}(A))=\lim_{n}\f_{\xi}(\g_{g_{n}}(A))
=\lim_{n}\f_{\xi}(A)\equiv\f_{\xi}(A)\,,
$$
which, by separability, implies 
$G_{\f}\subset G_{\f_{\xi}}$ whenever $\xi\in F_{0}$. Set
$F:=F_{0}\cap F_{1}$. Taking into account the first part of the proof 
and the definition \eqref{ddf} of $\f$, we obtain that $\xi\in F$ implies
$$
G_{\f_{\xi}}\subset G_{\f}\subset G_{\f_{\xi}}=G_{1}\,,
$$
which is the assertion.
\end{proof}
\begin{thm}
\label{gaga1}
If $\f\in\cs_{0}(\gf)$ is a weakly clustering state whose restriction to 
$\ga$ is $\gpt$--invariant, then there exist a continuous one--parameter 
subgroup $t\in\br\mapsto\eps_{t}\in\cz(G_{\f})$, and a measurable subset $F$ of 
full measure such that $\xi\in F$ implies that $\f_{\xi}$ is invariant 
under the modified time translation
$$
\f_{\xi}(A)=\f_{\xi}(\t^{\xi}_{t}\g_{\eps_{t}}(A)),
$$
simultaneously for every $A\in\cf$ and $t\in\br$.
\end{thm}
\begin{proof}
We cannot directly apply Theorem II.2 of \cite{AKTH} as our time 
translations $\gpt$ enjoy less continuity property than strong 
continuity, see \cite{BF}. However, in order to apply the mentioned result,
it is enough to verify that the one parameter group 
$t\in\br\mapsto\bar{v}_{t}\in\cn(G_{\f})/G_{\f})$ in pag. 106 of \cite{AKTH}
is continuous also in our situation (i.e. when the map 
$t\mapsto\f(A\gpt_{t}(B)C)$ is continuous for every fixed elements 
$A,B,C\in\gf$). Suppose not. There exists an open neighbourhood 
$U\supset G_{\f}$ of the stabilizer $G_{\f}$ of $\f$ such that 
$v_{1/n}\in U^{c}$. Choose a subsequence $\{v_{1/n_{k}}\}$ converging 
to some element $v_{0}\in U^{c}$. We get
$$
\f(A)=\lim_{k}\f(\gpt_{1/n_{k}}(A))=\lim_{k}\f(\gpg_{v_{1/n_{k}}}(A))
=\f(\gpg_{v_{0}}(A))
$$
which is a contradiction as the automorphism $\gpg_{v_{0}}$ is not
in the stabilizer. Hence, the conclusions of Theorem II.2 of 
\cite{AKTH} hold true also in our situation. We can conclude by 
reasoning as in the proof of Proposition \ref{gaga}, and after choosing 
countable dense subsets $\br_{0}\subset\br$, $\cf_{0}\subset\cf$. 
Namely, there exists a 
measurable subset $F\subset X$ of full measure such that $\xi\in F$, 
$t\in\br_{0}$ and $A\in\cf_{0}$ 
implies $\f_{\xi}(A)=\f_{\xi}(\t^{\xi}_{t}\eps_{t}(A))$. Let $t\in\br$ and 
$A\in\cf$. Choose convergent sequences $t_{n}\to t$, $A_{n}\to A$. We 
have on $F$, taking into account that 
$\t^{\xi}_{t_{n}}\eps_{t_{n}}(A_{n})\to \t^{\xi}_{t}\eps_{t}(A)$,  
$$
\f_{\xi}(A)=\lim_{n}\f_{\xi}(A_{n})
=\lim_{n}\f_{\xi}(\t^{\xi}_{t_{n}}\eps_{t_{n}}(A_{n}))
=\f_{\xi}(\t^{\xi}_{t}\eps_{t}(A))
$$
which is the assertion.
\end{proof}

Let $\f\in\cs_{0}(\gf)$, and consider the corresponding 
equivariant field $\{\f_{\xi}\}_{\xi\in X}\subset\cs(\cf)$. We have 
shown in Theorem \ref{stbb} that $G_{\f_{\xi}}=G_{\f}$ almost surely. 
This means that the GNS representations $\pi_{\f_{\xi}}$, as well as 
$\pi_{\f}$, are equipped with a strongly continuous representation
$U_{\f_{\xi}}$, or $U_{\f}$, of the common subgroup $G_{\f}\subset G$ 
implementing the gauge 
action of $G_{\f}$ on $\cf$ or $\gf$ respectively. 
Let $H\subset G_{\f}$ be a closed subgroup. Denote $\widehat{H}$ 
the set of all irreducible representations of $H$. It is 
well--known that the elements of $\widehat{H}$ act on finite 
dimensional Hilbert spaces (compactness of $H$), and 
$\widehat{H}$ is at most countable (second countability of $H$).

The restriction of
$U_{\f_{\xi}}$, or $U_{\f}$ to $H$ are denoted as 
$U^{H}_{\f_{\xi}}$, or $U^{H}_{\f}$ respectively. The $H$--spectra
$\S^{H}_{\f_{\xi}}, \S^{H}_{\f}\subset\widehat{H}$ of $\f_{\xi}$, $\f$ 
are the set of all 
irreducible representations of $H$ contained in
$U^{H}_{\f_{\xi}}$, or $U^{H}_{\f}$ respectively.

Following Definition II.3 of  \cite{AKTH}, we say that 
$\S^{H}_{\f_{\xi}}$ (or equivalently $\S^{H}_{\f}$) is {\it 
one--sided} if it is contained in a set $\S\subset\widehat{H}$ which 
enjoys the following properties:
\begin{itemize}
\item[(i)] $\s_{1},\s_{2}\in\S$ implies that every irreducible 
summand of $\s_{1}\otimes \s_{2}$ is also contained in $\S$, 
\item[(ii)] $\s,\bar\s\in\S$ implies that $\s=\id$.
\end{itemize}
\begin{thm}
\label{hsp}
Let $\f\in\cs_{0}(\gf)$. Then
\begin{itemize}
\item[(i)] the $H$--spectrum $\S^{H}_{\f_{\xi}}$ is almost surely 
independent on $\xi\in X$,
\item[(ii)] if $\s\in\S^{H}_{\f_{\xi}}$ almost surely, its 
multiplicity is (almost surely) independent on $\xi\in X$.
\end{itemize}
\end{thm}
\begin{proof}
We can identify the GNS triplet 
$(\pi_{\f_{T_{-x}\xi}},\ch_{\f_{T_{-x}\xi}},
\Om_{\f_{T_{-x}\xi}})$ relative to $\f_{T_{-x}\xi}$ with 
$(\pi_{\f_{\xi}}\circ\a_{x},\ch_{\f_{\xi}},\Om_{\f_{\xi}})$ almost 
surely. Under this identification, $U^{H}_{\f_{T_{-x}\xi}}$ coincides 
with $U^{H}_{\f_{\xi}}$. In other words, 
$U^{H}_{\f_{T_{-x}\xi}}\cong U^{H}_{\f_{\xi}}$ almost surely. 

The assertion follows by ergodicity, as the measurable subsets
$$
F_{\s,m}:=\big\{\xi\in X\,\big|\,\s\prec U^{H}_{\f_{\xi}}\, 
\text{with multiplicity}\,m\big\}
$$
give a countable partition of $X$.
\end{proof}
\begin{cor}
\label{hsp1}
Let $\f\in\cs_{0}(\gf)$.
\begin{itemize}
\item[(i)] If $\f$ is weakly clustering and asymptotically Abelian, 
then $U^{H}_{\f_{\xi}}$ 
satisfies the semigroup property of Theorem II.3 of \cite{AKTH} almost 
surely,
\item[(ii)] $\S^{H}_{\f}$ is one--sided if and only if $\S^{H}_{\f_{\xi}}$
is one--sided almost surely.
\end{itemize}
\end{cor}
\begin{proof}
The proof easily follows by Theorem \ref{hsp} and the mentioned 
Theorem II.3 of \cite{AKTH}.
\end{proof}

Here, there is the main theorem of the present section concerning the 
appearance of the chemical potential in the setting of disordered 
systems.
\begin{thm}
\label{gaga2}
If $\f\in\cs_{0}(\gf)$ is a weakly clustering asymptotically abelian
state whose restriction to 
$\ga$ is $(\gpt,\b)$--KMS state at inverse temperature $\b\neq0$, 
then there exist a closed subgroup $N\subset G_{\f}$, a continuous one--parameter 
subgroup $t\in\br\mapsto\eps_{t}\in\cz(G_{\f})$, a continuous one--parameter 
subgroup $t\in\br\mapsto\z_{t}\in G_{\f}$, and a measurable subset $F$ of 
full measure such that, for each $\xi\in F$, 
\begin{itemize}
\item[(i)] the $N$--spectrum of $\f_{\xi}$ is one--sided,
\item[(ii)] the restriction of $\f_{\xi}$ to 
$\cf^{N}:=\big\{A\in\cf\,\big|\,\g_{g}(A)=A\,,g\in N\big\}$ is a 
$(\th^{\xi},\b)$--KMS state for the modified time evolution 
$\th^{\xi}_{t}:=\t^{\xi}_{t}\g_{\eps_{t}\z_{t}}$,
\item[(iii)] the image $[\z_{t}]:=\z_{t}N$ in $G_{\f}/N$ is in $\cz(G_{\f}/N)$.
\end{itemize}
\end{thm}
\begin{proof}
We start by noticing that the conclusions of Theorem II.4 of \cite{AKTH} 
hold true also in our situation. Indeed, the hypotesis of extremality 
w.r.t the KMS condition is not used in the proof of that 
theorem.\footnote{See 
also the analogous results Theorem 4.1 of \cite{KT1}, and 
Theorem 12 of \cite{KT2}.}
Furthermore, in order to apply those results to our situation, we should 
replace the dense subset of entire elements used in II.6, pag. 110 
with the dense subset $\gf_{0}$ generated by elements of the form
$$
A_{f}(\xi):=\int f(t)\t^{\xi}_{t}(A(\xi))\di t\,,
$$
where $A$ runs on elements of $\gf$ and $f\in\widehat{\cd}$, 
together with (the image in $\gf$ of) $L^{\infty}(X,\n)$. 
Hence, we can apply the above mentioned theorem 
to $\f$ as above. The assertion follows with $N:=N_{\f}$,
$\eps$ and $\z$, the one parameter 
subgroups relative to $\f$ as in Theorem II.4 of \cite{AKTH}, 
by applying Theorem \ref{stbb}, Theorem 
\ref{gaga1}, Corollary \ref{hsp1}, and finally Proposition 3.2 of \cite{BF}.
\end{proof}

The subgroup $N_{\f}\subset G_{\f}$ of Theorem II.4 of \cite{AKTH} is 
called {\it the asymmetry subgroup} of $\f$.
Notice that, by an elementary application of Corollary \ref{hsp1}, 
in the situation of Theorem \ref{gaga2} the asymmetry subgroup 
$N_{\f_{\xi}}\subset G_{\f_{\xi}}$ of $\f_{\xi}$ coincides with  
the asymmetry subgroup $N_{\f}$ of $\f$ almost everywhere.

\section{an intrinsic characterization of the chemical potential for 
disordered systems}
\label{secint}

In the present section we show how the chemical potential can arise as an 
object directly associated to the algebra of observables. 
To simplify matter, we consider the simplest non trivial case when the gauge group is the 
unit circle $\bt$. In this 
situation, the charges in the model under 
consideration are generated by the powers $[\s^{n}]\in\text{Out}(\ca)$ of 
a single localized transportable automorphism $\s$, see \cite{DHR1}.\footnote{One can 
directly start from the observable algebra $\ca$, and then reconstruct the 
field algebra $\cf$ by considering the localizable charges of interest of 
the model, see \cite{DHR2, DR2}. This picture applies also to the case 
described in \cite{NS}, obtaining models satisfying all the properties 
assumed in the present section.}

We start by proving the following
\begin{prop}
\label{esdi}
Let $\om\in\cs_{0}(\ga)$ be strongly clustering asymptotically abelian $(\gpt,\b)$--KMS 
state at inverse temperature $\b\neq0$. If $\ca$ is simple, 
then any localized automorphism $\r$ of $\ca$ extends to a 
pointwise--weak measurable field $\{\r_{\xi}\}_{\xi\in X}$ of normal 
automorphism  of the weak closure $\pi_{\om_{\xi}}(\ca)''$, 
almost surely. 
\end{prop}
\begin{proof}
Under our assumptions, we can apply Proposition IV. 1 of \cite{AKTH}. 
Then $\om$ is quasi--equivalent to 
$\om\circ(\r\otimes\id)$, where 
$\id\equiv\id_{L^{\infty}(X,\n)}$.\footnote{In Proposition IV. 1 of \cite{AKTH},
the extremality of $\om$ cannot be dropped, as no weak clustering 
assumption is made there, see Remark II.2 of \cite{AKTH}.} Then 
$\pi_{\om}$ is unitarily equivalent to 
$\pi_{\om\circ(\r\otimes\id)}\equiv\pi_{\om}\circ(\r\otimes\id)$. 
This means that there exists a spatial isomorphism of
$$
\pi_{\om\circ(\r\otimes\id)}(\ga)''=\int^{}_{X}\pi_{\om_{\xi}}(\r(\ca))''
\n(\di\xi)
$$
onto
$$
\pi_{\om}(\ga)''=\int^{}_{X}\pi_{\om_{\xi}}(\ca)''\n(\di\xi)\,,
$$
where, under the above identification, we have a common
direct integral decomposition $\{\ch_{\xi}\}_{\xi\in X}$
of the same Hilbert space $\ch_{\om}$ on which $\pi_{\om}$ and 
$\pi_{\om\circ(\r\otimes\id)}$ act simultaneously.\footnote{See \cite{BF}, 
Section 2 for the last equality relative to the direct integral 
decomposition.} By Theorem IV.8.23 of \cite{T1},
There exists  a measurable field $\{U_{\xi}\}_{\xi\in X}$ of unitary 
operators such that
$$
\pi_{\om_{\xi}}(\r(\ca))''=U_{\xi}\pi_{\om_{\xi}}(\ca)''U_{\xi}^{*}\,,
$$
with
$$
\pi_{\om_{\xi}}\circ\r=U_{\xi}\pi_{\om_{\xi}}(\,\cdot\,)U_{\xi}^{*}
$$
almost everywhere. As $\ca$ is supposed to be simple, 
$\{\pi_{\om_{\xi}}\}_{\xi\in X}$ is a measurable field of 
$*$--isomorphism of $\ca$ onto their ranges $\pi_{\om_{\xi}}(\ca)$ 
almost surely.

After identifying $\ca$ with $\pi_{\om_{\xi}}(\ca)$, the measurable 
field of normal automorphisms are given, for $R\in\pi_{\om_{\xi}}(\ca)''$, by
$$
\r_{\xi}(R):=U_{\xi}RU_{\xi}^{*}\,.
$$
\end{proof}

In order to avoid technical problems, we suppose also that the quasi--local 
algebra of observables $\ca$ (or equally well the field algebra $\cf$  
is the $C^{*}$--inductive limit of local algebras isomorphic to a 
common full matrix algebra (i.e. the spin algebra). 

Let $\om\in\cs_{0}(\ga)$ be a $(\gpt,\b)$--KMS state 
such that the centre 
$\gz_{\pi_{\om}}:=\pi_{\om}(\ga)'\bigwedge\pi_{\om}(\ga)''$ is isomorphic 
to $L^{\infty}(X,\n)$. In \cite{BF} it is explained that this situation 
seems to be the right one in order to describe the ``pure 
termodynamical phase'' in the case of disordered models. 
In this situation, we have that $\om$ is weakly clustering (w.r.t. the 
spatial translation). Namely, as the $(\t^{\xi},\b)$--KMS state
$\om_{\xi}\in\cs(\ca)$ is a factor state almost surely, it satisfies
\eqref{ces} almost surely, see \cite{A}, Lemma 10.2.\footnote{See 
\cite{R}, Proposition 3, when $\ca$ includes Fermion operators, even 
if the last situation is not the standard one (e.g. \cite{DHR1, DR2}).}
Then the assertion follows by Proposition \ref{ces1}. 

Take a weakly 
clustering estension $\f$ of $\om$ to all of $\gf$ which  
exists by Proposition \ref{wce}. Furthermore, in order to avoid cases 
when the chemical potential is zero, we suppose that $\f$ is 
gauge--invariant (i.e. $G_{\f}=\bt$). In this situation, we have 
that $\sp(U_{\f})=\bz\equiv\widehat{G}$.\footnote{Let $V_{\s}\in\cf$ be the unitary 
implementing the automorphism $\s$ on $\ca$. This means that 
$\g_{\th}(V_{\s})=e^{i\th}V_{\s}$. The vector 
$\Psi_{n}:=\pi_{\f}(V_{\s}^{n}\otimes I)\Om_{\f}$ is an eigenvector 
of $U_{\f}(\th)$ corresponding to the eigenvalue $e^{in\th}$.} 
Namely, the asymmetry subgroup $N_{\f}$ of $\f$ is 
trivial.\footnote{By the 
previous results, one conclude that also $N_{\f_{\xi}}$ is trivial 
almost surely.}

Let $\r$ be a localized automorphism 
of $\ca$ carrying the charge $n$ (i.e. $\r\in[\s^{n}]$), and consider 
the unitary $U$ implementing $\r$ on $\ca$. Consider the state
$\f_{U}:=\f\circ\text{ad}_{U\otimes I}$. We have for the 
Connes--Radon--Nikodym cocycle (\cite{C, St}),
\begin{align}
\label{crncoc}
\big(D\f_{U}:D\f)_{t}=&\int^{\oplus}_{X}
\pi_{\f_{\xi}}(U^{*})\s^{\f_{\xi}}_{t}(\pi_{\f_{\xi}}(U))\n(\di\xi)\nn\\
=e^{in\b\m t}&\int^{\oplus}_{X}
\pi_{\f_{\xi}}(U^{*}\t^{\xi}_{-\b t}(U))\n(\di\xi)\,,
\end{align}
for some $\m\in\br$.

Here, we have used $\g_{\th}(U)=e^{in\th}U$, and
$\s^{\f}_{t}\circ\pi_{\f}=\pi_{\f}\circ\gpt_{-\b t}\circ\gpg_{\b\m t}$, 
by Theorem \ref{gaga2} taking into account \eqref{modgns}. 

Now, we take advantage from the fact that $\om\circ(\r\otimes\id)$ 
extends to a normal state on all of $\gf$. Denote with an abuse of 
notation, $\{\om_{\xi}\circ\r_{\xi}\}_{\xi\in X}$ the (equivariant) 
measurable field of states providing the direct integral decomposition of such an extension. 
Here, the $\r_{\xi}$ are the normal automorphisms of $\pi_{\f_{\xi}}$ 
appearing in Proposition \ref{esdi}. We have, by 
\eqref{crncoc} and the fact that $U^{*}\t^{\xi}_{-\b t}(U)$ is 
gauge--invariant,
\begin{equation}
\label{crncoc1}
\big(D(\om_{\xi}\circ\r_{\xi}):D\om_{\xi})=
e^{in\b\m t}
\pi_{\om_{\xi}}(U^{*}\t^{\xi}_{-\b t}(U))\
\end{equation}
almost everywhere.

Formula \eqref{crncoc1} explains the occurrence of the 
chemical potential $\m\in\br$ as an object intrinically associated to the 
observable algebra. Furthermore, according with this description, it
does not depend on the disorder for states 
$\om$ on $\ga$ such that $\gz_{\pi_{\om}}\sim L^{\infty}(X,\n)$.
This is 
in accordance with standard
fact that the physically relevant quantities should not depend on
the disorder. 

As it is explained in \cite{AKTH}, Section IV, there is a freeness 
in order to define the chemical potential (see Formula (IV.6)). 
By this freeness, the chemical potential might be defined up to a 
phase factor in the centre of the GNS 
representation of the state. Such a phase is connected with the zero--point 
of the chemical potential, as that centre is trivial in the 
situation treated in the above mentioned paper. In our situation, 
such a zero--point would lie in $L^{\infty}(X,\n)$. However, if one make 
such a choice in a measurable and 
invariant way, one would conclude by ergodicity, that also the 
freeness in the choice of the zero--point of the chemical potential 
can be avoided.

To conclude, the following remarks are in order. First, the situation described in the 
present section extends 
straighforwardly to the case when the gauge group is the 
$n$--dimensional torus. Second, one could extend the matter to more 
complicated 
situations arising from continuous disordered systems,\footnote{See 
\cite{FL1, FL2} for some non trivial examples of continuous 
disordered models.} as well as 
possible disordered systems arising from quantum field theory. In 
these cases, one could take advantage from the local normality of the 
objects of interest and/or the split property naturally assumed in 
quantum field theory.\footnote{See \cite{L1}, Section 3 for an 
analysis of the chemical potential for some (non disordered) models of 
low dimensional quantum field theory, without passing from the field 
algebra.} We choose not to pursue such an 
analysis as, at the knowledge of the author, no natural 
disordered model arising from quantum field theory seems to be 
present in literature.

\end{document}